\documentclass{article}
\usepackage{amsmath}
\usepackage{amssymb}

\newtheorem{theorem}{Theorem}[section]
\newtheorem{lemma}[theorem]{Lemma}

\newcommand{\QED}{\hfill$\Box$\\~}

\makeatletter
 
  \@addtoreset{equation}{section}
\makeatother

\begin{document}
\title{\bf Generalized eigenfunctions of relativistic Schr\"odinger operators  \\
in two dimensions}
\author{Tomio Umeda\thanks{Supported by Grant-in-Aid for Scientific Research (C) No.09640212 and (A) No.19204013.
from the Japan Society for the Promotion of Science.}  \ and Dabi Wei}

\date{}
\maketitle
\begin{abstract}



Generalized eigenfunctions of the two-dimensional 
relativistic Schr\"o\-dinger operator 
$H=\sqrt{-\Delta}+V(x)$ with $|V(x)|\leq C\langle x\rangle^{-\sigma}$, 
$\sigma>3/2$, are considered. 
We compute the integral kernels of the boundary values 
$R_0^\pm(\lambda)=(\sqrt{-\Delta}-(\lambda\pm i0))^{-1}$, 
and prove that 
the generalized eigenfunctions $\varphi^\pm(x,k)$ 
are bounded on $R_x^2\times\{k\,| \,a\leq |k|\leq b\}$, 
where $[a,b]\subset(0,\infty)\backslash\sigma_p(H)$, and 
$\sigma_p(H)$ is the set of eigenvalues of $H$. 
With this fact and the completeness of the wave operators,
we establish the eigenfunction expansion for the absolutely continuous subspace
for $H$.
Finally, we show that 
each generalized eigenfunction is asymptotically equal to
 a sum
of a plane wave and a spherical wave
 under the assumption that $\sigma>2$.
\end{abstract}

\section{Introduction}   \label{sec:intro}

Generalized eigenfunctions for 
 Schr\"odinger
operators $-\Delta+V(x)$ on ${\mathbb R}^n$ are  now
well understood at least in the framework
of simple scattering; see for example Agmon\cite{Agmon1},
Ikebe\cite{Ikebe2} and Kato and Kuroda\cite{KatoKuroda}.
In the pseudo-relativistic regime, 
one can replace the Schr\"odinger operators
with relativistic 
Schr\"odinger operators
$\sqrt{-\Delta+m}+V(x)$.
Here $m$ is the mass of the particle, and it could be zero.
In this case, we deal with
the operators
of the form $\sqrt{-\Delta}+V(x)$. 

This paper is a continuation of our 
previous paper Wei \cite{Wei1}, where
 the odd-dimensional relativistic Schr\"odinger operators $\sqrt{-\Delta}+V(x)$
were considered
and  substantial generalizations of the results by Umeda \cite{Umeda5}, \cite{Umeda3},
who only dealt with the three-dimensional case, 
were accomplished.
In the present paper, we shall deal with 
 the two-dimensional case: 
\begin{equation}
H=H_0+V(x), \quad H_0=\sqrt{-\Delta}, 
\quad x\in \mathbb R^2.
\label{eq:defH}
\end{equation}
Our aim  here is to establish all the same 
results as in \cite{Umeda5}, \cite{Umeda3} and \cite{Wei1}. 
For this reason and for the sake of simplicity, we shall 
use the same notation as in \cite{Wei1}.

We now roughly recall the discussions demonstrated in
our previous works \cite{Umeda5}, \cite{Umeda3} and \cite{Wei1}
for the reader's convenience.
We first defined the generalized eigenfunctions  $\varphi^{\pm}(x,k)$ 
with the aid of the limiting absorption principle
for the relativistic Schr\"odinger operators.
 We next proved that the generalized eigenfuctions
are bounded 
on the set $(x,k)\in\mathbb R^n\times\{k\in\mathbb R^n \,| \,a\leq |k|\leq b\}$ 
for $[a,b]\subset(0,\infty)\backslash\sigma_p(\sqrt{-\Delta}+V(x))$, 
where $n= 3$, $5$, $7$, $\cdots$, and $\sigma_p(\sqrt{-\Delta}+V(x))$
denotes the point spectrum. 
Then we showed the asymptotic completeness of the wave operators
by the Enss method (cf. \cite{Enss1,Isozaki1}), 
and obtained the eigenfunction expansions for the absolutely continuous subspace 
for $\sqrt{-\Delta}+V(x)$. 
In the three dimensional case,
 we gave estimates on the differences between the
generalized eigenfunctions and the plane waves.  Moreover we showed that 
the generalized eigenfunctions are
asymptotically equal to the sum of  plane waves and  spherical waves.
It should be  remarked that
once we have the boundedness of the generalized eigenfunctions, 
we are able to establish the completeness of the generalized eigenfunctions
for the absolutely continuous subspace. 
(See \cite{Wei1}. Also see \cite{Kitada3, Kuroda1}.)

\vspace{5pt}

Our basic assumption is as follows.

\vspace{7pt}

\noindent
\textbf{Assumption:}
$V(x)$ is a
real-valued measurable  function on $\mathbb R^2$ satisfying
\begin{equation}  
|V(x)|\leq C\langle x\rangle^{-\sigma}, \quad
\sigma>3/2.
\label{eq:assumV}
\end{equation}

\vspace{5pt}

Under the assumption (\ref{eq:assumV}), it is obvious that  
$V= V(x) \times$ is a bounded selfadjoint
operator in $L^2(\mathbb R^2)$, and
that $H=H_0+V$ defines a selfadjoint operator in
$L^2(\mathbb R^2)$,  whose domain is $H^1(\mathbb R^2)$,
the Sobolev space of order one. 
  Moreover $H$ is essentially
selfadjoint on $C_0^\infty(\mathbb R^2)$ (see \cite[sections 2 and 7\,]{Umeda3}).  
Note that
\[
\sigma_{e}(H)=\sigma_{e}(H_0),
\]
where $\sigma_{e}(H)$ and $\sigma_{e}(H_0)$ denote
the essential spectrum of $H$ and $H_0$ respectively.
This fact follows from Reed and Simon \cite[p.113, Corollary
2]{Reed1},
since $V$ is
relatively compact with respect to $H_0$. 
Also, note that the essential spectrum of $H_0$ coincides with 
the spectrum of $H_0$: $\sigma_{e}(H_0)=\sigma(H_0)=[0, \, +\infty)$.

The main idea in this paper
 is essentially the same as in \cite{Umeda3}
and \cite{Wei1}. Thus we basically follow the same
line as in \cite{Umeda3} and \cite{Wei1}.
Namely,  we first prove the boundedness of the generalized eigenfunctions, 
and then we  establish the eigenfunction expansion,
and finally we examine asymptotic behaviors of 
the generalized eigenfunctions at infinity.

However, we should like to emphasize that 
some difficulties specific to the two-dimensional case
arise.
One should recall  that 
there are significant differences between 
the two-dimensional wave equation and 
the three-dimensional one in their 
treatments.
We find that
a similar phenomenon is also observed in the treatments
of relativistic 
Schr\"odinger
operators.

In the odd-dimensional case, the integral kernel of 
the resolvent of the operator $\sqrt{-\Delta}$ 
is expressed in terms of 
trigonometric functions, and the cosine and sin 
integral functions (see \cite{Umeda3} and \cite{Wei1}). 
On the other hand,
we encounter 
the Bessel function,
Neumann function and the Struve function
in the integral kernel of the resolvent of  $\sqrt{-\Delta}$ 
in the two-dimensional case.
This difference makes the analysis of
the resolvent of  $\sqrt{-\Delta}$ 
 more difficult in the two-dimensional case.
 
In fact, when we deal with the boundary values of the
the resolvent $R_0(z)$ 
to define the generalized eigenfunctions of
$H$,  we are obliged to
examine the boundary values of 
 all of the Bessel function,
the Neumann function and the Struve function
on the positive half line $[0, \, +\infty)$.
It is surprising that
a suitable  combination of these special functions
on the positive half line exhibits a 
simple form of an exponential function at infinity.
This fact enables us to show  that
the generalized eigenfunctions of 
relativistic Schr\"odinger operators in the two-dimensional case too
are asymptotically equal to superpositions of plane waves
and spherical waves at infinity .

We would like to mention a technicality.
In showing the boundedness of generalized eigenfunctions 
in section \ref{sec:boundedness},
we need to handle the Riesz potential on ${\mathbb R}^2$.
We shall show that for functions in a certain class 
the Riesz potential defines bounded functions. 
We believe that this fact, as well as our technique,
 is interesting in its own right.
The key for this fact is the estimate (\ref{eq:etvarphi7}),
which is based on  Lemma
\ref{lm:D10psiInLQ}.

We expect that the discussions on the generalized eigenfunctions in 
the $2m$ dimensional case ($m \ge 2$)  would become more
complicated, and will be discussed elsewhere (\cite{Wei3}).

\vspace{10pt}

\textbf{The plan of the paper}~~
In section \ref{sec:GE}, we define generalized eigenfunctions of $H$.
In section \ref{sec:kernel},
  we compute the resolvent kernel of $H_0$. 
Section \ref{sec:boundedness} is devoted to prove
the boundedness of the generalized eigenfunctions. 
In  section \ref{sec:expansions}, 
we deal with the completeness of the generalized eigenfunctions
for the absolutely continuous subspace for $H$. 
Finally, in section \ref{sec:asymptotic}, 
we examine the asymptotic behaviors of the generalized eigenfunctions
at infinity.
In appendix, we include two inequalities which are used repeatedly
in the present paper, and summarize some basic properties
of the Bessel, Neumann and Struve functions for the
reader's convenience. 

\section{Generalized eigenfuctions}  \label{sec:GE}

By $R(z)$ and $R_0(z)$, we mean the resolvents of $H$ and 
$H_0$ respectively:
\begin{equation}
R(z):=(H-z)^{-1}, \;\; R_0(z):=(H_0-z)^{-1}.
\end{equation}

The task of this section is to
construct  generalized eigenfunctions $\varphi^{\pm}(x, \, k)$  of $\sqrt{-\Delta}+V(x)$
(see Theorem \ref{th:defPhi1} below), 
and  show that they satisfy
\begin{equation}  \label{eqn:UdefGE1}
\varphi^{\pm}(x,k)=\varphi_0(x,k)-R_0^{\mp}(|k|)V\varphi^{\pm}(x,k),
\end{equation} 
where $R_0^{\pm}(z)$ denotes the extended resolvents of
$H_0$ (cf. Theorem \ref{th:1.1} below) 
and $\varphi_0(x,k)$ denotes the plane wave
\begin{equation}
\varphi_0(x,k)=e^{ix\cdot k}.
\end{equation}
As we shall see in Theorem \ref{Th:kernel} 
in section \ref{sec:kernel},
the extended resolvents $R_0^{\pm}(\lambda)$ have 
the integral kernels $g_\lambda^\pm(x-y)$.
Since we have 
\begin{equation*}
g_\lambda^\pm(x)
\approx
\left(\frac\lambda\pi\right)^{1/2}(1\mp i)
\frac{e^{\mp i(\lambda|x|-\pi/4)}}{|x|^{1/2}}
\end{equation*}
as $|x|\to \infty$ (see (\ref{eq:g_|x|2infty}) below),
 it is justified to call (\ref{eqn:UdefGE1})  the Lippman-Schwinger
type integral equations.

The discussions in this section are based
on the results  by
 Ben-Artzi and Nemirovski \cite[sections 2 and 4]{Ben1}.
Since their results are formulated in a general setting,
we reproduce them in the context of 
the present paper.

\begin{theorem}[Ben-Artzi and Nemirovski \cite{Ben1}]
Let $s>1/2$. Then\\
\mbox{\rm(1)}~~For any $\lambda>0$, there exist the limits 
$R_0^{\pm}(\lambda)=\lim_{\mu\downarrow0}R_0(\lambda\pm i\mu)$ 
in ${\textbf B}(L^{2,s},H^{1,-s})$.\\
\mbox{\rm(2)} The operator-valued functions $R_0^{\pm}(z)$ defined by
\[
R_0^\pm(z)=\begin{cases}R_0(z)~~~\mbox{if}~~~z\in\mathbb 
C^\pm\\R_0^\pm(\lambda)~~~\mbox{if}~~~z=\lambda>0\end{cases}
\]
are ${\textbf B}(L^{2,s},H^{1,-s})$-valued continuous functions, 
where $\mathbb C^+$ and $\mathbb C^-$ are the upper and
 the lower half-planes respectively: 
$\mathbb C^\pm=\{z\in\mathbb C \, |\, \pm \mbox{\rm Im }z>0\}$. 
\label{th:1.1}
\end{theorem}

\begin{theorem}[Ben-Artzi and Nemirovski \cite{Ben1}]
Let $s>1/2$ and $\sigma>1$. Then\\
\mbox{\rm(1)}~~The continuous spectrum $\sigma_c(H)=[0,\infty)$ is 
absolutely continuous, 
except possibly for a discrete set of embedded 
eigenvalues $\sigma_p(H)\cap(0,\infty)$, 
which can accumulate only at $0$ and $\infty$.\\
\mbox{\rm(2)}~~For any $\lambda\in(0,\infty)\backslash\sigma_p(H)$, 
there exist the limits
\[
R^{\pm}(\lambda)=\lim_{\mu\downarrow0}R(\lambda\pm i\mu)~~~~
\mbox{in}~~~~{\textbf B}(L^{2,s},H^{1,-s}).
\]
\mbox{(3)} The functions $R^{\pm}(z)$ defined by
\[
R^\pm(z)=\begin{cases}R(z)~~~\mbox{if}~~~z\in\mathbb C^\pm\\
R^\pm(\lambda)~~~\mbox{if}~~~z=\lambda\in(0,\infty)\backslash\sigma_p(H)\end{cases}
\]
are ${\textbf B}(L^{2,s},H^{1,-s})$-valued continuous.
\label{th:1.2}\end{theorem}

Now we can follow the arguments in our previous papers
 \cite[section 8]{Umeda3} and \cite[section 1]{Wei1} 
with a few of obvious changes, and 
obtain the following two theorems.

\begin{theorem}[\cite{Umeda3}, \cite{Wei1}]
If $|k|\in(0,\infty)\backslash\sigma_p(H)$, 
then the eigenfunctions defined by
\begin{equation}
\varphi^\pm(x,k)=\varphi_0(x,y)-R^{\mp}(|k|)\{V(\cdot)\varphi_0(\cdot,k)\}(x),
\label{eq:defvarphi1}
\end{equation}
satisfy the equation
\[(\sqrt{-\Delta_x}+V(x))u=|k|u ~~~in~~~\mathcal S'(\mathbb R_x^2).\]
\label{th:defPhi1}
\end{theorem}

\begin{theorem}[\cite{Umeda3}, \cite{Wei1}]
If $|k|\in(0,\infty)\backslash\sigma_p(H)$ 
and $1<s<\sigma-1/2$, then we have
\[\varphi^\pm(x,k)=\varphi_0(x,k)-R_0^\mp(|k|)\{V(\cdot)\varphi^\pm(\cdot,k)\}(x)~~~
in~~~L^{2,-s}(\mathbb R^2).\]
\label{th:defPhi2}
\end{theorem}

\section{The integral kernels of  the resolvents of $H_0$}   \label{sec:kernel}

This section is devoted to computing the kernel 
$g_z(x-y)$
of the resolvent $R_0(z)$. 
What we shall need in the later sections is the limit  
$g_{\lambda}^{\pm}(x)$ of the function $g_{\lambda\pm i\mu}(x)$ as
$\mu\downarrow0$,  where $\lambda>0$. 
Then we derive a few inequalities for the extended resolvent
$R_0^{\pm}(\lambda)$,  using  some estimates of the functions
 $g_{\lambda}^{\pm}(x)$.

We  first need to introduce the 
following functions.
\begin{align}    
 M_z(x)&=\frac1{2}\Big\{{\mathbf H}_0(-|x|z)-N_0(-|x|z)\Big\},
  \quad 
  z \in {\mathbb C} \setminus  [0,\, +\infty) \label{eqn:LY-1}\\
\noalign{\vskip 2pt}
m_\lambda^\pm(x)&= -\frac12\Big\{ {\mathbf H}_0(|x|\lambda)+N_0(|x|\lambda) 
\pm 2iJ_0(|x|\lambda)\Big\}, 
\quad  \lambda >0.   \label{eqn:LY-2}
\end{align}
Here ${\mathbf H}_0(z)$ 
is the Struve function (cf. \cite[p.227, p.228]{Moriguti2}, \cite[p.328]{Watson1}),
$N_0(z)$  the Neumann function (cf.  \cite[p.145, p.146]{Moriguti2}, 
\cite[p.62,
p.64]{Watson1}; the Neumann function is  denoted by $Y_0(z)$  
in \cite{Watson1}) and
 $J_0(z)$  the Bessel function  (cf. \cite[p.145,
p.146]{Moriguti2}, \cite[p.40]{Watson1}):
\begin{align}
{\mathbf H}_0(z)&=\sum_{k=0}^\infty
\frac{(-1)^k(z/2)^{2k+1}}{\{\Gamma(k+3/2)\}^2},\label{eq:H0}\\
N_0(z)&=\frac2\pi J_0(z)(\gamma+\log(z/2))-
\frac2\pi\sum_{k=1}^\infty\frac{(-1)^k(z/2)^{2k}\sum_{m=1}^k\frac1m}{(k!)^2}
\label{eq:N0} \\
\noalign{\vskip -3pt}
 & \;\; (\gamma \;\; \mbox{  the Euler constant}).   \nonumber \\
\noalign{\vskip 3pt}
J_0(z)&=\sum_{n=0}^\infty\frac{(-1)^n(z/2)^{2n}}{(n!)^2}.\label{eq:J0}
\end{align}
Note that the Struve function ${\mathbf H}_0(z)$ and
the Bessel function $J_0(z)$ are both entire functions.
Also note that the Neumann function $N_0(z)$ is a many-valued function with 
a logarithmic branch-point at $z=0$. 
Here we choose the principal branch, i.e. 
$|\mbox{Im }\log z| <\pi$ 
for $z \in {\mathbb C}\setminus (-\infty, \, 0]$.

The resolvent kernel of $H_0$ is given as follows.
\begin{theorem}
If $z \in {\mathbb C}\setminus [0, \, +\infty)$, then
\[
R_0(z)u=G_zu
\]
for all $u\in C_0^\infty(\mathbb R^2)$, where
\begin{align}
G_zu(x)&=\int_{\mathbb R^2}g_z(x-y)u(y)dy,   \\
\noalign{\vskip 4pt}
g_z(x)&=\frac1{\pi|x|}+zM_z(x).  \label{eqn:Uresolvent37-3}
\end{align}
\label{th:resolvent1}
\end{theorem}

\noindent
{\it Proof.} We follow the same line as in
\cite[senction 2]{Umeda3} and \cite[section 2]{Wei1},
and we only give the sketch of the proof.

We start with the Poison kernel 
\begin{equation*}
P_t(x) = \frac{t}{\pi(t^2+|x|^2)^{3/2}},
\end{equation*}
and the fact that $e^{-tH_0}u = P_t *u$ for $t>0$ and
 $u\in L^2({\mathbb R}^2)$. Then we appeal to the fact that
\begin{equation*}
R_0(z)=\int_0^{+\infty}e^{tz}e^{-tH_0}dt, \quad
\mbox{Re } z <0.
\end{equation*}
For all $u$, $v\in C_0^\infty(\mathbb R^2)$ we have
\begin{align}
&( R_0(z)u, \, v)_{L^2}    \nonumber\\
& = 
  \int_0^{+\infty} e^{tz} \, 
( e^{-tH_0}u, \, v)_{L^2} \, dt  \nonumber\\  
 &= 
 \,\int_{{\mathbb R}^2}   \!
        \Big\{ \int_{{\mathbb R}^2}  \!
           \Big( \int_0^{+\infty} e^{tz}  
\frac{t}{\pi(t^2+|x-y|^2)^{3/2}} \, dt  \Big) \,
u(y) \, dy \Big\}  
\overline{v(x)} \, dx 
\label{eqn:U2y10}                
\end{align}
for $z$ with $\hbox{Re }z <0$, where
 we have made a change of order of integration.
 (Note that the integral in (\ref{eqn:U2y10}) is absolutely
 convergent. See the proof of \cite[Lemma 2.2]{Wei1}, 
 which is valid in any dimension $n \ge 2$.)
It is evident that
 the integration with respect $t$ in (\ref{eqn:U2y10}) gives
the integral kernel of $R_0(z)$ if $\hbox{Re }z <0$.
For this reason we make the following computation:
\begin{equation}
\begin{split}
&\int_0^\infty e^{tz}\frac{t}{\pi(t^2+|x|^2)^{3/2}}dt\\
&=\left[-\frac{e^{tz}}{\pi\sqrt{t^2+|x|^2}}\right]_0^\infty+z
\int_0^\infty 
\frac{e^{tz}}{\pi\sqrt{t^2+|x|^2}}dt\\
&=\frac1{\pi|x|}+zM_z(x)  \\
&=g_z(x)
\end{split}
\end{equation}
if $\hbox{Re }z <0$.
Here we have used the formula 
\[
\int_0^\infty \frac{e^{tz}}{\pi\sqrt{t^2+|x|^2}}dt
=\frac1{2}\Big\{ {\mathbf H}_0(-|x|z)-N_0(-|x|z)\Big\}
\]
(cf. \cite[p.138]{Erdelyi}, \cite[p.289]{Moriguti1};
note that the Neumann function is  denoted by $Y_0(z)$  
in \cite{Erdelyi}).
Summing up, we have shown that
\begin{equation}
( R_0(z)u, \, v)_{L^2} = ( G_z u, \, v)_{L^2}
\label{eqn:U39resolvent}
\end{equation}
for all $u$, $v\in C_0^\infty(\mathbb R^2)$
when $\hbox{Re }z <0$.
Since both sides of (\ref{eqn:U39resolvent}) are
holomorphic functions of $z$ on ${\mathbb C}\setminus [0, \, +\infty]$,
we get the conclusion of the theorem.
\QED

In the proof of Theorem \ref{Th:kernel} below,
we shall need the following estimates (Appendix B):
For $\rho>0$
\[\begin{split}
&|J_0(\rho)|\leq \text{const.}\begin{cases}1 
     & \mbox{if } \; 0<\rho\leq 1, \\ \rho^{-1/2} & \mbox{if } \; \rho\geq 1,\end{cases}\\
\noalign{\vskip 4pt}
&|N_0(\rho)|\leq \text{const.}\begin{cases} |\log \rho| 
        &\mbox{if } \;  0<\rho\leq 1, \\ \rho^{-1/2} 
        & \mbox{if } \;  \rho\geq 1,\end{cases}\\
\noalign{\vskip 4pt}
&| {\mathbf H}_0(\rho)|\leq \text{const.}\begin{cases} \rho 
         & \mbox{if } \; 0< \rho \leq 1, \\ \rho^{-1/2}
         & \mbox{if } \; \rho \geq 1 .\end{cases}
\end{split}\]
Since
 $|\log \rho| \leq \text{const.} \rho^{-1/2}~~~(0<\rho\leq 1)$, 
we see that 
\begin{equation}
|m_\lambda^\pm(x) |\leq \text{const.} (|x|\lambda)^{-1/2}.
\label{eqn:Umlambda}
\end{equation}

\vspace{5pt}

\begin{theorem}
If $\lambda >0$, then
\[
R_0^\pm(\lambda)u=G_\lambda^\pm u
\]
for all $u\in C_0^\infty(\mathbb R^2)$, where
\begin{align}
\begin{split}
G_\lambda^\pm u(x)
&=\int_{\mathbb R^2}g_{\lambda}^{\pm}(x-y)u(y)dy,   \\
\noalign{\vskip 4pt}
g_{\lambda}^{\pm}(x)&=\frac1{\pi|x|}+\lambda m_\lambda^\pm(x).
\end{split}
\label{eq:R0_kernel}
\end{align}
\label{Th:kernel}
\end{theorem}

\vspace{5pt}
\noindent
{\it Proof.}  Again we follow the same line as in
\cite[senction 4]{Umeda3} and \cite[section 2]{Wei1},
and we only give the sketch of the proof.

It follows from Theorem \ref{th:resolvent1} that
\begin{equation}
( R_0(\lambda \pm i \mu)u, \, v)_{L^2} 
 = ( G_{\lambda \pm i \mu} u, \, v)_{L^2}
\label{eqn:U39resolvent-1}
\end{equation}
for all $u$, $v\in C_0^\infty(\mathbb R^2)$
whenever $\lambda>0$, $\mu >0$.
Regarding  $R_0(\lambda \pm i \mu)u \in L^{2, -s}$
 and $v \in L^{2, s}$ for some $s>1/2$,
we apply  Theorem \ref{th:1.1} to the left-hand
side of  (\ref{eqn:U39resolvent-1}), and see that
\begin{equation}
\lim_{\mu \downarrow 0} ( R_0(\lambda \pm i \mu)u, \, v)_{L^2} 
=( R_0^{\pm}u, \, v)_{-s, s}.
\label{eqn:U39resolvent-2}
\end{equation}
Here $( \cdot, \, \cdot)_{-s, s}$ denotes the anti-duality bracket or 
the pairing between
$L^{2, -s}$ and $L^{2, s}$.
To examine the limit of the right-hand side
of (\ref{eqn:U39resolvent-1}),
we see that
\begin{align*}
\begin{split}
&\lim_{\mu \downarrow 0} {\mathbf H}_0(-|x|(\lambda\pm i\mu)) 
    =-{\mathbf H}_0(|x|\lambda),   \\
&\lim_{\mu \downarrow 0} J_0(-|x|(\lambda\pm i\mu)) =J_0(|x|\lambda),   \\
&\lim_{\mu \downarrow 0} 
N_0(-|x|(\lambda\pm i\mu))
 =N_0(|x|\lambda)\pm 2i J_0(|x|\lambda).
\end{split}
\end{align*}
These facts, together with (\ref{eqn:LY-1}), (\ref{eqn:LY-2}), 
(\ref{eqn:Uresolvent37-3})
and  (\ref{eq:R0_kernel}), 
show that
\begin{equation}
\lim_{\mu\downarrow0}g_{\lambda\pm i\mu}(x)
=\frac1{\pi|x|}+\lambda m_\lambda^\pm(x)
=g_\lambda^\pm(x).
\end{equation}
By virtue of (\ref{eqn:Umlambda}),
 we can apply the Lebesgue dominated convergence theorem
to the right-hand side
of (\ref{eqn:U39resolvent-1}), and
we get
\begin{equation}
\lim_{\mu\downarrow 0}( G_{\lambda \pm i \mu} u, \, v)_{L^2}
=
\iint_{\!\!{\mathbb R}^4}
g_{\lambda}^{\pm} (x-y) u(y) \overline{v(y)} \, dx \, dy.
\label{eqn:URO-kernel-316}
\end{equation}
Combining (\ref{eqn:U39resolvent-2}) and (\ref{eqn:URO-kernel-316}),
we get the conclusion of the theorem.
\QED

\vspace{5pt}

It follows from Theorem \ref{Th:kernel} that
the integral operator $G_{\lambda}^{\pm}$ can be
extended to bounded operators
from $L^{2, s}$ to $H^{1, -s}$ for $s > 1/2$.

We shall show the boundedness of the generalized
eigenfunctions $\varphi^{\pm}(x,k)$ 
in section \ref{sec:boundedness}, where we shall use 
the following integral operators: 
\begin{equation}
T_ju(x):=\int_{\mathbb R^2}|x-y|^{-j}u(y)dy, \quad j=1, \,1/2
\label{eqn:UT12}
\end{equation}
Recall that these integral operators are 
actually Riesz potentials up to constants.

Following lemma is a direct consequence of Theorem \ref{Th:kernel},
(\ref{eqn:Umlambda})and (\ref{eqn:UT12}).

\vspace{5pt}
\begin{lemma} Let $s > 1/2$.
If $[a, \, b] \subset (0, \, +\infty)$.
then there exist a positive  constant $C_{ab}$ such that
\begin{equation}
|R_0^\pm(\lambda) u(x)|
\leq 
\frac{1}{\pi} \big| T_1u  (x) \big|
  +C_{ab} \big( T_{1/2}|u| \big) (x)
\label{lem:UT123}
\end{equation}
for all $u \in L^{2,s}$ and all $\lambda \in [a, \, b]$.
\label{lem:lemma31}
\end{lemma}

We prepare one more lemma for a later purpose.

\begin{lemma}
For each $\lambda >0$ we have
\begin{equation}
g_\lambda^\pm(x)
=\left(\frac\lambda\pi\right)^{1/2}(1\mp i)
\frac{e^{\mp i\lambda|x|}}{|x|^{1/2}}+O(|x|^{-1})
\label{eq:g_|x|2infty}
\end{equation}
as $|x|\to\infty$.
\label{lem:lemma32}
\end{lemma}

\vspace{5pt}
\noindent
\textit{Proof.}
Apply
Lemmas \ref{lm:J0N0inf} and \ref{lm:H0inf} in the appendix
to (\ref{eqn:LY-2}).
\QED

\section{Boundedness of the generalized eigenfunctions}   
\label{sec:boundedness}

In this section, we shall discuss 
the boundedness of the 
generalized eigenfuctions $\varphi^{\pm}(x,k)$ defined in Theorem \ref{th:defPhi1}.
Following our previous papers \cite{Umeda3} and \cite{Wei1}, 
we shall need a restriction on $k$.
Namely, we 
assume that $k$ satisfies the following inequality:
\begin{equation}
a\leq|k|\leq b,
\label{eq:assumk}
\end{equation}
 where $[a,b]\subset(0,\infty)\backslash\sigma_p(H)$.
As we have seen in Theorem \ref{th:defPhi2},  the
generalized eigenfuction $\varphi^{\pm}(x,k)$
satisfies the equation
\begin{equation}
\varphi^\pm(x,k)
=\varphi_0(x,k)-R_0^\mp(|k|)\{V(\cdot)\varphi^\pm(\cdot,k)\}(x).
\label{eq:defvarphi2}
\end{equation}
In section \ref{sec:kernel}, we have shown that $R_0^\mp(|k|)$ 
are integral operators, and investigated properties
of the integral kernels.

We are now in a position to state the main theorem in this section,
which is stated as follows.

\begin{theorem}
Let $[a,b]\subset(0,\infty)\backslash\sigma_p(H)$. 
There exists a constant $C_{ab}$ such that 
generalized eigenfunctions defined by (\ref{eq:defvarphi1}) satisfy
\begin{equation}
|\varphi^\pm(x,k)|\leq C_{ab}
\label{eq:varphiBounded}\end{equation}
for all $(x,k)\in\mathbb R^2\times\{\, k \, | \, a\leq|k|\leq b\,\}$.
\label{th:varphiBond}
\end{theorem}

Before proving Theorem \ref{th:varphiBond},
we have to prepare a few lemmas.
With application of Theorem \ref{Th:kernel} in mind, 
we shall show that $V(x)\varphi^{\pm}(x,k)$ belongs to 
$L^{2,s}(\mathbb R_x^2)$  
provided that $1/2<s < \sigma -1$. 
To this end, we put
\begin{equation}
\psi^\pm(x,k)=V(x)\varphi^\pm(x,k).
\label{eqn:umevarphi44}
\end{equation}

\begin{lemma}
If $1/2<s < \sigma -1$, 
then $\psi^\pm(x,k)$ are $L^{2,s}(\mathbb R_x^2)$--valued continuous functions on 
$\big\{ \, k \, \big| \, |k|\in(0,\infty)\backslash\sigma_p(H) \big\}$.
\label{lm:psiInL2s}
\end{lemma}

Proof.~~Since we have \cite[Lemma 1.1]{Wei1} with $n=2$,
we can imitate the arguments in
\cite[Lemmas 9.2 and 9.3]{Umeda3}, and  see that
for any $t>1$,
$\psi^\pm(x,k)$ are $L^{2,\sigma-t}(\mathbb R_x^2)$-valued
continuous functions 
on 
$\big\{ \, k \, \big| \, |k|\in(0,\infty)\backslash\sigma_p(H) \big\}$.
For $s \in (1/2, \, \sigma -1)$, we put $t:= \sigma -s$.
Then $t >1$, and hence we get the lemma.
\QED

\begin{lemma}
If  $4/3 < r < 2$, then 
$\psi^\pm(x,k)$ are $L^{r}(\mathbb R_x^2)$--valued continuous fuctions on 
$\big\{ \, k  \, \big| \, |k|\in(0,\infty)\backslash\sigma_p(H) \big\}$.
\label{lm:psiInL16/9Ume}
\end{lemma}

Proof.~~Applying the H\"older inequality, we have
\[
\begin{split}
\int_{\mathbb R^2}|\psi^\pm(x,k)|^{r}dx
\leq&
\left\{\int_{\mathbb R^2}
\left(\langle x\rangle^{-r/2}\right)^{2/(2-r)}dx\right\}^{(2-r)/2}\\
&\times\left\{\int_{\mathbb R^2}
\left(\langle x\rangle^{r/2}|
\psi^\pm(x,k)|^{r}\right)^{2/r}dx\right\}^{r/2}\\
=&
\left\{\int_{\mathbb R^2}
\langle x\rangle^{-r/(2-r)}dx\right\}^{(2-r)/2}
\left\{\int_{\mathbb R^2}
\langle x\rangle
\, |\psi^\pm(x,k)|^2 dx\right\}^{r/2}\\
=&
C_r \left(\|\psi^\pm\|_{L^{2,1/2}}\right)^{r/2} \\
\leq&
 C_r\left(\|\psi^\pm\|_{L^{2,s}}\right)^{r/2} < \infty,
\end{split}
\]
where $C_r$ is a constant depending only on $r$ 
and $s\in (1/2, \, \sigma -1)$. 
Here we have used the fact that $r/(2-r) > 2$ if 
 and only if $4/3 < r <2$.
Lemma \ref{lm:psiInL2s}, together with this inequality,
implies that 
$\psi^\pm(x,k)$ belongs to $L^{r}(\mathbb R_x^2)$ 
if  $4/3 < r < 2$. 
Moreover, 
by using a similar argument,
one can easily show that 
$\psi^\pm(x,k)$ are $L^{r}(\mathbb R_x^2)$-valued 
continuous functions 
on $\big\{ \, k  \, \big| \, |k|\in(0,\infty)\backslash\sigma_p(H) \big\}$.
\QED

For the sake of simplicity,
we shall apply Lemma \ref{lm:psiInL16/9Ume} with $r=16/9$: 
\begin{lemma}
$\psi^\pm(x,k)$ are $L^{16/9}(\mathbb R_x^2)$--valued continuous functions on 
$\big\{ \, k  \, \big| \, |k|\in(0,\infty)\backslash\sigma_p(H) \big\}$.
\label{lm:psiInL16/9}
\end{lemma}

\vspace{5pt}

As we mentioned in section \ref{sec:kernel}, 
we shall use the integral operators $T_1$ 
and $T_{1/2}$; see (\ref{eqn:UT12}).
It will be convenient 
to  split $T_1$ into two parts:
\begin{equation}
T_1=T_{10}+T_{1\infty},
\label{eq:defDj-Ume}
\end{equation}
where
\begin{eqnarray*}
T_{10}u(x)&=&\int_{|x-y|\leq1}|x-y|^{-1}u(y) \, dy,\\
T_{1\infty}u(x)&=&\int_{|x-y|>1}|x-y|^{-1}u(y) \, dy.
\end{eqnarray*}
Then it follows from Lemma \ref{lem:lemma31} 
 and (\ref{eq:defvarphi2}) that
\begin{equation}
\begin{split}
{}&|\varphi^\pm(x,k)|
\leq  1 +  C_{ab} \\
& \times 
\Big\{
\big| \big( T_{10}\psi^\pm(\cdot,\,k) \big) (x) \big|
+ \big| \big( T_{1\infty}\psi^\pm(\cdot,\,k) \big) (x) \big| 
 + \big( T_{1/2}
|\psi^\pm(\cdot, \, k)|\big)(x)
\Big\}
\label{eq:etvarphi1}
\end{split}
\end{equation}
for all $(x, \, k) \in 
{\mathbb R}^2 \times
\big\{\, k \, | \,a\le |k| \le b \big\}$,
where $C_{ab}$ is a positive constant.

\vspace{3pt}

\begin{lemma}
If $16/9\leq q<16$, then 
$T_{10}\psi^\pm(\cdot, \,k)\in L^{q}(\mathbb R^2)$. 
Moreover, there exits a positive constant $C_{ab}$ such
that
\[
\Vert 
T_{10}\psi^\pm(\cdot, \,k)
\Vert_{L^q}
\leq C_{ab}
\]
for all   $k \in \big\{ \, k \, | \,a\le |k| \le b \big\}$.
\label{lm:D10psiInLQ}
\end{lemma}

Proof.~~ We write
\[
\big( T_{10}\psi^\pm(\cdot,k) \big)(x)
=\int_{\mathbb R^2}f_0(x-y)\psi^\pm(y,k)dy,~~
f_0(x):=|x|^{-1}\chi_0(x),
\]
where $\chi_0(x)$ is the characteristic
function for the unit disk $\{\, x \, | \, |x|\leq 1\}$.
It is easy to see that 
\begin{equation}
f_0\in L^p(\mathbb R_x^2) \text{ for all } p \in (0, \, 2). 
\label{eq:f0InLP}
\end{equation}
Using Lemma \ref{lm:psiInL16/9} and 
the Young inequality (cf. Lemma
\ref{lm:YoungInequality} in the appendix) with $r=16/9$,  
we get 
\[
\Vert 
T_{10}\psi^\pm(\cdot,k)
\Vert_{L^q}
\leq\|f_0\|_{L^p}\|\psi^\pm(\cdot,k)\|_{L^{16/9}}
\]
for $\frac1q=\frac1p+\frac9{16}-1$ $(1\leq p,q\leq\infty)$.
Noticing (\ref{eq:f0InLP}), we have
\[
\frac1{16}<\frac1q\leq\frac9{16}
\Longleftrightarrow 
\frac{16}9\leq q<16.
\]
Thus we get the lemma.
\QED

We are now in a position to prove the main theorem  in this section,
 namely Theorem
\ref{th:varphiBond}.
In the proof  below,
we shall apply Lemma \ref{lm:D10psiInLQ} 
with  $q=3$.

\vspace{10pt}
{\bf Proof of Theorem \ref{th:varphiBond}}
Let $1/2<s<\sigma-1$.
Noticing the definition (\ref{eqn:UT12}) and the Schwarz inequality,
we have 
\[
\big( T_{1/2}|\psi^\pm(\cdot,k)|\big)(x)
\leq
\left\{
\int_{\mathbb
R^2}\frac1{|x-y|\langle y\rangle^{2s}}dy
\right\}^{1/2}
\left\{
\int_{\mathbb R^2}\langle y\rangle^{2s}|\psi^\pm(y,k)|^2dy
\right\}^{1/2}.
\]
Using Lemma \ref{lm:psiInL2s} and Lemma \ref{lm:UmedaA.1}
in the appendix with
$\beta=1,~\gamma=2s>1,~n=2$,  we get
\begin{equation}
|\big( T_{1/2}|\psi^\pm(\cdot,k)|\big)(x)|
\leq C_{ab1}'
\label{eq:etvarphi2}
\end{equation}
for all $(x, \, k) \in 
{\mathbb R}^2 \times
\big\{\, k \, | \,a\le |k| \le b \big\}$,
where $C_{ab1}'$ is a positive constant. 

Lemma \ref{lm:psiInL16/9},
together with the H\"older inequality, yields
\begin{equation}
\begin{split}
|T_{1\infty}\psi^\pm(x,k)|
\leq&
\left\{  \int_{|x-y|>1}|x-y|^{-16/7} \, dy \right\}^{7/16}  \\
&\times
\left\{ \int_{\mathbb R^2}|\psi^\pm(y,k)|^{16/9}\, dy \right\}^{9/16} \\
\leq&C_{ab2}'
\end{split}
\label{eq:etvarphi3}
\end{equation}
for all $(x, \, k) \in 
{\mathbb R}^2 \times
\big\{\, k \, | \,a\le |k| \le b \big\}$,
where $C_{ab2}'$ is a positive constant.

Combining  (\ref{eq:etvarphi2}), 
(\ref{eq:etvarphi3}) and (\ref{eq:etvarphi1}), we have 
thus shown that
\begin{equation}
\begin{split}
|\varphi^{\pm}(x, \, k)|
\leq&
1 + C_{ab} 
\big\{ C_{ab1}'+C_{ab2}' + \big(T_{10}|\psi^\pm(\cdot,k)|\big)(x) \big\} \\
\noalign{\vskip 4pt}
=&
C_{ab}'' 
\big\{
1 + 
\big(T_{10}|V(\cdot)\varphi^\pm(\cdot,k)|\big)(x) 
\big\}.
\end{split}
\label{eq:etvarphi3UME}
\end{equation}
(Recall (\ref{eqn:umevarphi44}).)
Here we would like to utilize the fact
that $T_{10}$ is positivity preserving, i.e.
\begin{equation}
T_{10}u \ge 0 \mbox{\  if } u \ge 0.
\label{eqn:positivity}
\end{equation}
It then follows from (\ref{eq:etvarphi3UME}) and  (\ref{eqn:positivity})
that
\begin{equation}
\begin{split}
|\varphi^{\pm}(x, \, k)|
\leq&
C_{ab}'' 
\Big\{
1 + 
\Big( T_{10}
|V(\cdot)|
C_{ab}'' 
\big\{
1 + 
\big(T_{10}|V(\cdot)\varphi^\pm(\cdot,k)|\big) 
\big\}
\Big)(x) 
\Big\}  \\
\noalign{\vskip 4pt}
=
C_{ab}'' 
\Big\{
1 
&+ 
C_{ab}'' 
\big( T_{10}|V(\cdot)| \big)(x)    \\
&\quad+ 
C_{ab}'' 
\Big( T_{10}|V(\cdot)| 
\big( T_{10}|\psi^\pm(\cdot,k)|\big) 
\Big)(x) 
\Big\} 
\end{split}
\label{eq:etvarphi3UME+1}
\end{equation}
(Again recall (\ref{eqn:umevarphi44}).)

With the same notation as in the proof of
Lemma \ref{lm:D10psiInLQ}, we have
\begin{equation}
0 \leq \big( T_{10}|V(\cdot)| \big)(x)  
\leq \Vert f_0 \Vert_{L^{3/2}} \Vert V \Vert_{L^3} < +\infty,
\label{eqn:V(x)+ume}
\end{equation}
where we have used the H\"older inequality.

Similarly, by using the H\"older inequality and 
applying Lemma \ref{lm:D10psiInLQ},
we have
\begin{equation}
\begin{split}
0 \le& 
T_{10}|V(\cdot)| 
\big( T_{10}|\psi^\pm(\cdot,k)|\big) 
\Big)(x)     \\
\leq&
\Vert f_0 \Vert_{L^{3/2}} 
\Vert V \Vert_{L^\infty}
\Vert T_{10}|\psi^\pm(\cdot,k)| \Vert_{L^3}   \\
\leq&C_{ab3}'
\end{split}
\label{eq:etvarphi7}
\end{equation}
for all $(x, \, k) \in 
{\mathbb R}^2 \times
\big\{\, k \, | \,a\le |k| \le b \big\}$,
where $C_{ab3}'$ is a positive constant.

Combining (\ref{eq:etvarphi3UME+1}) with (\ref{eqn:V(x)+ume}) and (\ref{eq:etvarphi7}),
we obtain the desired conclusion.
\QED


\section{Generalized eigenfunction expansions} \label{sec:expansions}

The task in this section is
to establish the completeness of the generalized eigenfunction.  
The idea is
the same as  in our previous work \cite{Wei1}. 
For this reason, we shall only state the results
and  omit the proofs.

It is obvious that $V$ is a bounded
selfadjoint operator in $L^2(\mathbb R^2)$,  
and that $H=H_0+V$ defines a
selfadjoint operator in $L^2(\mathbb R^2)$,  
whose domain is $H^1(\mathbb R^2)$
(see \cite[Theorem 5.8]{Umeda2}).  
Moreover $H$ is essentially selfadjoint on
$C_0^\infty(\mathbb R^2)$  (see  \cite{Umeda2}). 
Since $V$ is relatively compact with respect to $H_0$, 
it follows from \cite[p.113, Corollary 2]{Reed1} that
\[
\sigma_{e}(H)=\sigma_e(H_0)=[0,\infty).
\]

The first result in this section is the
asymptotic completeness of wave operators
 (cf. \cite{Wei1}).

\begin{theorem}
Let $H_0,H$ be defined by \mbox{\rm(\ref{eq:defH})} and $V(x)$ satisfy 
\mbox{\rm(\ref{eq:assumV})}. 
Then there exist the limits
\[
W_{\pm}=\mbox{\rm s-}\!\!\!\!
\lim_{t\to\pm\infty}e^{itH}e^{-itH_0},
\]
and the asymptotic completeness holds:
\[
\mathcal R(W_{\pm})=\mathcal H_{ac}(H),
\]
\label{th:AsymptoticCompleteness}
where $\mathcal H_{ac}(H)$ denotes the 
absolutely continuous subspace for $H$.
\end{theorem}

We need to remark that $\sigma_p(H)\cap(0,\infty)$ is a discrete set. 
This fact was first proved by B. Simon \cite[Theorem 2.1]{Simon1}. 
Moreover, he proved that each eigenvalue in the set 
$\sigma_p(H)\cap(0,\infty)$ has finite multiplicity.
Finally, using Theorem \ref{th:varphiBond} and Theorem
\ref{th:AsymptoticCompleteness},  
we can establish the eigenfunction expansion theorem as follows
(see our previous work \cite{Wei1} for the details).

\begin{theorem}
Let $H_0,H$ be defined by \mbox{\rm(\ref{eq:defH})} and 
$V(x)$ satisfy \mbox{\rm(\ref{eq:assumV})}. 
Let $s>1$ and $[a,b]\subset(0,\infty)\backslash\sigma_p(H)$. 
For $u\in L^{2,s}(\mathbb R^2)$, let $\mathcal F_{\pm}$ be defined by
\[
\mathcal F_{\pm}u(k):=(2\pi)^{-1} 
\int_{\mathbb
R^2}u(x)\overline{\varphi^\pm(x,k)}dx.
\]
Then for any   $f  \in L^{2,s}(\mathbb R^2)$, 
we have
\[
E_H([a,b])f(x)
=(2\pi)^{-1}\int_{a\leq|k|\leq b}\mathcal F_\pm f(k)
\varphi^\pm(x,k)dk,
\] 
where
$E_H$ is the spectral measure of $H$.
\label{th:mainTh}
\end{theorem}


\section{Asymptotic behaviors of the generalized eigenfunctions}   \label{sec:asymptotic}

We shall first show that the generalized eigenfunctions $\varphi^{\pm}(x,k)$, 
defined by (\ref{eq:defvarphi1}), 
are distorted plane waves, 
and give estimates of the differences between $\varphi^{\pm}(x,k)$ 
and the plane wave $\varphi_0(x,k)=e^{ix\cdot k}$ 
(Theorem \ref{th:asymp2planwave}).
We shall next prove that $\varphi^{\pm}(x,k)$ are asymptotically equal to 
the sums of the plane wave and the spherical waves $e^{\mp i|x||k|}/|x|^{1/2}$ 
under the assumption that $\sigma>2$, 
and shall give estimates of the differences between $\varphi^{\pm}(x,k)$ and 
the sums mentioned above (Theorem \ref{th:asymp2planwave+sphwave}). 

The similar estimates were discussed in T. Ikebe \cite[\S 3]{Ikebe2} and 
our previous work
\cite[\S 10]{Umeda3},  though our arguments below are slightly different from those
of \cite{Ikebe2} or \cite{Umeda3},  and our estimates are slight refinements of
those of \cite{Ikebe2} or \cite{Umeda3}.

The main theorems in this section are

\begin{theorem}
Let $\sigma>3/2$. If $|k|\in(0,+\infty)\backslash\sigma_p(H)$, then
\[
|\varphi^{\pm}(x,k)-\varphi_0(x,k)|\leq C_{k}
\begin{cases}
\langle x\rangle^{-(\sigma-3/2)}& \mbox{ if}~~ 3/2<\sigma<2,\\
\langle x\rangle^{-1/2}\log(1+\langle x\rangle)& \mbox{ if}~~ \sigma=2,\\
\langle x\rangle^{-1/2}& \mbox{ if}~~ \sigma>2.
\end{cases}
\]
where the constant $C_k$ is uniform for $k$ in any compact subset of 
\[\left\{k\Big||k|\in(0,+\infty)\backslash\sigma_p(H)\right\}.\] 
\label{th:asymp2planwave}\end{theorem}

\begin{theorem}
Let $\sigma>2$ and 
\begin{equation}
f^\pm(\lambda,\omega_x,\omega_k)
:=
\left(\frac\lambda\pi\right)^{1/2}(1\mp i)
\int_{\mathbb R^2}
e^{\pm i\lambda \omega_x\cdot y}
V(y)\varphi^\pm(y,\lambda\omega_k)dy,
\label{eq:def_f_pm}\end{equation}
where $\omega_x=x/|x|$, $\omega_k=k/|k|$.
Then for $|x|\ge 1$
\begin{eqnarray}
&&\left|\varphi^\pm(x,k)-
\Big(
\varphi_0(x,k)+
\frac{e^{\mp i|k||x|}}{|x|^{1/2}}f^\pm(|k|,\omega_x,\omega_k)
\Big)
\right|    \nonumber  \\
\noalign{\vskip 4pt}
&&\quad \leq 
C_k
\begin{cases}
|x|^{-(\sigma-1)/2}& \mbox{if \ }2<\sigma <3, \\
|x|^{-1}& \mbox{if \ }\sigma \ge 3,
\end{cases} 
\label{eq:planwave+sphwave_esti}
\end{eqnarray}
where the constant $C_k$ is uniform for $k$ in any compact subset of 
\[
\left\{k\Big||k|\in(0,+\infty)\backslash\sigma_p(H)\right\}.
\] 
\label{th:asymp2planwave+sphwave}
\end{theorem}

We should like to remark that what makes the discussions below possible is 
the estimate in Theorem \ref{th:varphiBond}.

\vspace{5pt}
\noindent
{\bf Proof of Theorem \ref{th:asymp2planwave}}.~~
In view of (\ref{eq:defvarphi2}), Theorem \ref{Th:kernel} and 
Theorem
\ref{th:varphiBond}, it is clear that 
there is a positive constant $C_k$, which
is uniform for $k$ in any compact subset of 
$\{ \, k \, | \, |k|\in (0,+\infty)\backslash\sigma_p(H) \,\}$, 
such
that 
\begin{equation}
\begin{split}
|\varphi^{\pm}(x,k)-\varphi_0(x,k)|
\leq& 
C_k (T_1|V|)(x)+(T_{1/2}|V|)(x)  \\
\leq& 
C_k\Big(
\int_{\mathbb R^2}\frac1{|x-y|\langle y\rangle^{\sigma}}dy
+
\int_{\mathbb R^2}\frac1{|x-y|^{1/2}
\langle y \rangle^{\sigma}}dy
\Big).
\end{split}
\label{eq:pf_th:asymp2planwave_01}
\end{equation}
(Recall that $T_1$ and $T_{1/2}$ were introduced in
 (\ref{eqn:UT12}).)
We apply Lemma \ref{lm:UmedaA.1} 
with $n=2,\beta=1,\gamma=\sigma>3/2$, and get
\begin{equation}
\int_{\mathbb R^2}\frac1{|x-y|\langle y\rangle^{\sigma}}dy
\leq C_{\sigma}
\begin{cases}
\langle x\rangle^{-(\sigma-1)}& \mbox{\rm if}~~ 3/2<\sigma<2,\\
\langle x\rangle^{-1}\log(1+\langle x\rangle)& \mbox{\rm if}~~ \sigma=2,\\
\langle x\rangle^{-1}& \mbox{\rm if}~~ \sigma>2,
\end{cases}
\label{eq:pf_th:asymp2planwave_02}
\end{equation}
where $C_\sigma$ is a constant depending only on $\sigma$.
Similarly, 
we apply Lemma \ref{lm:UmedaA.1} 
with $n=2,\beta=1/2,\gamma=\sigma>3/2$,
and get
\begin{equation}
\int_{\mathbb R^2}\frac1{|x-y|^{1/2}\langle y\rangle^{\sigma}}dy
\leq C_{\sigma}
\begin{cases}
\langle x \rangle^{-(\sigma-3/2)}& \mbox{\rm if}~~ 3/2<\sigma<2,\\
\langle x \rangle^{-1/2}\log(1+\langle x\rangle)& \mbox{\rm if}~~ \sigma=2,\\
\langle x \rangle^{-1/2}& \mbox{\rm if}~~ \sigma>2.
\end{cases}
\label{eq:pf_th:asymp2planwave_03}
\end{equation}
The theorem is 
a direct consequence of (\ref{eq:pf_th:asymp2planwave_01}), 
(\ref{eq:pf_th:asymp2planwave_02}) and 
(\ref{eq:pf_th:asymp2planwave_03}).
\QED

We shall give a proof of Theorem \ref{th:asymp2planwave+sphwave} by means
of a series of lemmas. 

\vspace{5pt}
\begin{lemma}
Let $\sigma>2$. If
\begin{equation}
0\leq f(x) \leq C|x|^{-1},
\end{equation}
then
\begin{equation}
\int_{\mathbb R^2}f(x-y)\langle y\rangle^{-\sigma}dy=O(|x|^{-1})
\end{equation}
as $|x|\to \infty$, 
where $C$ is a constant.
\label{lm:et_f_-1/2in_-1out}\end{lemma}

Proof.~~ 
Applying Lemma \ref{lm:UmedaA.1} with $n=2,\beta=1,\gamma=\sigma>2$, 
we have
\begin{equation}\begin{split}
\int_{\mathbb R^2}f(x-y)\langle y\rangle^{-\sigma}dy
\leq &C\int_{\mathbb R^2}\frac1{|x-y|\langle y\rangle^{\sigma}}dy\\
\leq &C'\langle x\rangle^{-1},
\end{split}\label{eq:et_f_-1_out}\end{equation}
where $C$ and $C'$ are constants.
It is apparent that (\ref{eq:et_f_-1_out}) 
gives the lemma.
\QED

In view of (\ref{eqn:Umlambda}), (\ref{eq:R0_kernel}) and (\ref{eq:g_|x|2infty}),
we get 
\begin{equation}
\left|
g_\lambda^\pm(x)-\left(\frac\lambda\pi\right)^{1/2}(1\mp i)\frac{e^{\mp i\lambda|x|}}{|x|^{1/2}}
\right|
\leq C|x|^{-1},
\end{equation}
where $C$ is a constant.
Then, using Lemma \ref{lm:et_f_-1/2in_-1out}, 
(\ref{eq:defvarphi2}) and (\ref{eq:R0_kernel}),
we see that
\begin{eqnarray}
{}&&\varphi^\pm(x,k)-\varphi_0(x,k)  \nonumber \\
\noalign{\vskip 4pt}
&=&
\left(\frac\lambda\pi\right)^{\!\!1/2} \!\!
(1\mp i)
\int_{\mathbb R^2}
\frac{e^{\mp i|k| |x-y|}}{|x-y|^{1/2}}
V(y)\varphi^\pm(y,k)\, dy + O(|x|^{-1})
\label{eq:varphi-varphi0}
\end{eqnarray}
as $|x|\to\infty$.
Now, noticing  (\ref{eq:def_f_pm}), 
(\ref{eq:planwave+sphwave_esti}) and 
(\ref{eq:varphi-varphi0}), 
we need to consider the integral of 
the form
\begin{equation}
\int_{\mathbb R^2}
\left\{
\frac{e^{ia|x-y|}}{|x-y|^{1/2}}
-\frac{e^{ia(|x|-\omega_x\cdot y)}}{|x|^{1/2}}
\right\}
u(y)dy,
\label{eq:AB_keylemma_eq0}
\end{equation}
where $a\in\mathbb R$ and $u$ is a function satisfying
\begin{equation}
|u(x)|\leq C\langle x\rangle^{-\sigma},~~~\sigma>2.
\label{eqn:UMintu(x)}
\end{equation}

\vspace{5pt}

\begin{lemma}
Let  $u$ satisfy (\ref{eqn:UMintu(x)}).
Then for $|x|\geq1$ 
we have
\begin{eqnarray}
\left|\int_{|y|\geq\sqrt{|x|}}\frac{e^{ia(|x|-\omega_x\cdot y)}}{|x|^{1/2}}u(y)dy\right|
\leq C_1\|\langle\cdot\rangle^\sigma u\|_{L^\infty}|x|^{-(\sigma-1)/2},
\label{eq:AB_keylemma_eq1}\\
\left|\int_{|y|\geq\sqrt{|x|}}\frac{e^{ia|x-y|}}{|x-y|^{1/2}}u(y)dy\right|
\leq C_2\|\langle\cdot\rangle^\sigma u\|_{L^\infty}|x|^{-(\sigma-1)/2}.
\label{eq:AB_keylemma_eq2}
\end{eqnarray}
\label{lm:AB_keylemma_lm1}
\end{lemma}

Proof.~~
We obtain
\begin{eqnarray}
\left|\int_{|y|\geq\sqrt{|x|}}e^{ia(|x|-\omega_x\cdot y)}u(y)dy\right|
\leq C_1\|\langle\cdot\rangle^\sigma u\|_{L^\infty}|x|^{-(\sigma-2)/2}
\end{eqnarray}
by  similar arguments in \cite[(10.15)]{Umeda3}.
This inequality implies (\ref{eq:AB_keylemma_eq1}).

To prove (\ref{eq:AB_keylemma_eq2}), we write
\begin{eqnarray}
F_0(x):=\left\{ y\in \mathbb R^2 \left| \,
|y|\geq\sqrt{|x|},\,
|x-y|\leq \frac{|x|}2\right.\right\},  \\ 
F_1(x):=\left\{y\in \mathbb R^2 \left| \,
|y|\geq\sqrt{|x|}, \,
|x-y|\geq \frac{|x|}2\right.\right\},
\end{eqnarray}
and get
\begin{eqnarray}
\left|\int_{F_0}\frac{e^{ia|x-y|}}{|x-y|^{1/2}}u(y)dy\right|
\leq C'\|\langle\cdot\rangle^\sigma u\|_{L^\infty}|x|^{-(\sigma-3/2)},
\label{eq:AB_keylemma_eq3}  \\
\noalign{\vskip 4pt}
\left|\int_{F_1}\frac{e^{ia|x-y|}}{|x-y|^{1/2}}u(y)dy\right|
\leq C''\|\langle\cdot\rangle^\sigma u\|_{L^\infty}|x|^{-(\sigma-1)/2},
\label{eq:AB_keylemma_eq4}
\end{eqnarray}
by  similar arguments in 
\cite[(10.17) and (10.18)]{Umeda3}.

Since $\sigma>2\Leftrightarrow\sigma-3/2>(\sigma-1)/2$, 
we conclude from
(\ref{eq:AB_keylemma_eq3}) and (\ref{eq:AB_keylemma_eq4}) 
that the inequality (\ref{eq:AB_keylemma_eq2}) holds.
\QED

In view of (\ref{eq:AB_keylemma_eq0}) and 
Lemma \ref{lm:AB_keylemma_lm1}, 
it is sufficient to evaluate the integral
of the form
\begin{equation}
\int_{|y|\leq\sqrt{|x|}}\left\{\frac{e^{ia|x-y|}}{|x-y|^{1/2}}-
\frac{e^{ia(|x|-\omega_x\cdot y)}}{|x|^{1/2}}\right\}u(y)dy.
\end{equation}
We split it into two parts:
\begin{eqnarray}
&&\frac1{|x|^{1/2}}\int_{|y|\leq\sqrt{|x|}}\left\{e^{ia|x-y|}-e^{ia(|x|-\omega_x\cdot y)}\right\}
u(y)dy\nonumber\\
&&+\int_{|y|\leq\sqrt{|x|}}e^{ia|x-y|}\left(\frac1{|x-y|^{1/2}}-\frac1{|x|^{1/2}}\right)u(y)dy.
\end{eqnarray}
and evaluate these two integrals separately.

\begin{lemma}
If $\sqrt{|x|}\geq5$ and $|y|\leq \sqrt{|x|}$, then
\begin{equation}
\Big||x-y|-(|x|-\omega_x\cdot y)\Big|\leq 3\sqrt 2\frac{|y|^2}{|x|}.
\end{equation}
\label{lm:AB_keylemma_lm2}
\end{lemma}

For the proof of  this lemma, see \cite[(10.26)]{Umeda3}.

\begin{lemma}
Under the same assumptions as in Lemma \ref{lm:AB_keylemma_lm1}, 
we have
\begin{eqnarray}
&&\left|\frac1{|x|^{1/2}}\int_{|y|\leq\sqrt{|x|}}\left\{e^{ia|x-y|}-
e^{ia(|x|-\omega_x\cdot y)}\right\}u(y)dy\right|\nonumber\\
&&\leq C_3|a|~\|\langle\cdot\rangle^\sigma u\|_{L^\infty}
\begin{cases}
\displaystyle |x|^{-(\sigma-1)/2}&\mbox{if \ }2<\sigma<4,\\
\displaystyle |x|^{-3/2}\log(1+|x|)&\mbox{if \ }\sigma=4,\\
\displaystyle |x|^{-3/2}&\mbox{if \ }\sigma>4.
\end{cases}
\end{eqnarray}
for $\sqrt{|x|}\geq5$.
\label{lm:AB_keylemma_lm3}
\end{lemma}

Proof.~~
Let $\sqrt{|x|}\geq5$. 
In a similar fashion to
 in \cite[(10.28) and (10.30)]{Umeda3}, 
we get
\begin{eqnarray}
&&\left|\frac1{|x|^{1/2}}\int_{|y|\leq\sqrt{|x|}}
\left\{
e^{ia|x-y|}-e^{ia(|x|
-\omega_x\cdot y)}
\right\} u(y)dy
\right|           \nonumber\\
\noalign{\vskip 4pt}
&&\leq3\sqrt2|a|~
\|\langle\cdot\rangle^\sigma u\|_{L^\infty}\frac1{|x|^{3/2}}
\left|\int_{|y|\leq\sqrt{|x|}}|y|^2\langle y\rangle^{-\sigma}dy\right|
\label{eq:AB_keylemma_eq5}
\end{eqnarray}
and
\begin{eqnarray}
\left|\int_{|y|\leq\sqrt{|x|}}|y|^2\langle y\rangle^{-\sigma}dy\right|
&\leq&2^{\sigma/2}\int_0^{\sqrt{|x|}}(1+r)^{-\sigma+3}dr\nonumber\\
&\leq&
\begin{cases}
\displaystyle
2^{\sigma/2}\frac{|x|^{-(\sigma-4)/2}}{4-\sigma}& \mbox{if \ }2<\sigma<4,\\
\displaystyle 2^{\sigma/2}\log(1+|x|)& \mbox{if \ }\sigma=4,\\
\displaystyle 2^{\sigma/2}\frac1{4-\sigma}&\mbox{if \ }\sigma>4.
\end{cases}
\label{eq:AB_keylemma_eq6}
\end{eqnarray}
Combining
(\ref{eq:AB_keylemma_eq5}) with (\ref{eq:AB_keylemma_eq6}) yields the
desired inequalities.
\QED

\begin{lemma}
Under the same assumptions as in Lemma \ref{lm:AB_keylemma_lm1}, 
we have
\begin{eqnarray}
&&\left|\int_{|y|\leq\sqrt{|x|}}e^{ia|x-y|}
\left(\frac1{|x-y|^{1/2}}-\frac1{|x|^{1/2}}\right)u(y)dy\right|\nonumber\\
&\leq& C \, \|\langle\cdot\rangle^\sigma u\|_{L^\infty}
\begin{cases}
\displaystyle |x|^{-\sigma/2}& \mbox{if \ }2<\sigma<3,\\
\displaystyle |x|^{-3/2}\log(1+|x|)&\mbox{if \ }\sigma=3,\\
\displaystyle |x|^{-3/2}&\mbox{if \ }\sigma>3.
\end{cases}.
\end{eqnarray}
for $\sqrt{|x|}\geq5$.
\label{lm:AB_keylemma_lm4}
\end{lemma}

Proof.~~
It is follows that
\begin{eqnarray}
\left|\frac1{|x|^{1/2}}-\frac1{|x-y|^{1/2}}\right|
&=&\frac{\left||x-y|^{1/2}-|x|^{1/2}\right|}{|x|^{1/2}|x-y|^{1/2}}\nonumber\\
&=&\frac{\left||x-y|-|x|\right|}{|x|^{1/2}|x-y|^{1/2}\left||x-y|^{1/2}+|x|^{1/2}\right|}.
\label{eq:AB_keylemma_eq7}
\end{eqnarray}
If $\sqrt{|x|}\geq5$ and $|y|\leq\sqrt{|x|}$, 
then Lemma \ref{lm:AB_keylemma_lm2} implies
\begin{equation}
\big||x-y|-|x|\big|\leq|y|+3\sqrt2\frac{|y|^2}{|x|}.
\label{eq:AB_keylemma_eq8}
\end{equation}
If $\sqrt{|x|}\geq5$ and $|y|\leq\sqrt{|x|}$, 
 we then have
\begin{equation}
|x-y|\geq|x|-|y|\geq|x|-\frac{|x|}5=\frac45|x|.
\label{eq:AB_keylemma_eq9}
\end{equation}
Hence, 
it follows from (\ref{eq:AB_keylemma_eq7}), (\ref{eq:AB_keylemma_eq8}) 
and (\ref{eq:AB_keylemma_eq9}) 
that
\begin{eqnarray}
\left|\frac1{|x|^{1/2}}-\frac1{|x-y|^{1/2}}\right|
\leq C'\frac{|y|}{|x|^{3/2}}+C''\frac{|y|^2}{|x|^{5/2}}.
\label{eq:AB_keylemma_eq10}
\end{eqnarray}
when $\sqrt{|x|}\geq5$ and $|y|\leq\sqrt{|x|}$.
Using this inequality, we arrive at
\begin{eqnarray}
&&\left|\int_{|y|\leq\sqrt{|x|}}e^{ia|x-y|}
\left(\frac1{|x-y|^{1/2}}-\frac1{|x|^{1/2}}\right)u(y)dy\right|\nonumber\\
&\leq& C'\|\langle\cdot\rangle^\sigma u\|_{L^\infty}\frac1{|x|^{3/2}}
\int_{|y|\leq\sqrt{|x|}}|y|\langle y\rangle^{-\sigma}dy\nonumber\\
&&+C''\|\langle\cdot\rangle^\sigma u\|_{L^\infty}\frac1{|x|^{5/2}}
\int_{|y|\leq\sqrt{|x|}}|y|^2\langle y\rangle^{-\sigma}dy.
\label{eq:AB_keylemma_eq11}
\end{eqnarray}
provided that $\sqrt{|x|}\geq5$.
Also we have 
\begin{eqnarray}
\left|\int_{|y|\leq\sqrt{|x|}}|y|\langle y\rangle^{-\sigma}dy\right|
\leq\begin{cases}
\displaystyle 2^{\sigma/2}\frac{|x|^{-(\sigma-3)/2}}{3-\sigma}&2<\sigma<3,\\
\displaystyle 2^{\sigma/2}\log(1+|x|)&\sigma=3, \\
\displaystyle 2^{\sigma/2}\frac1{3-\sigma}&\sigma>3.
\end{cases}
\label{eq:AB_keylemma_eq12}
\end{eqnarray}
Combining (\ref{eq:AB_keylemma_eq11}) with (\ref{eq:AB_keylemma_eq12}) and 
(\ref{eq:AB_keylemma_eq6}), 
we conclude that the desired inequalities are verified.
\QED

\vspace{5pt}
\noindent
{\bf Proof of Theorem \ref{th:asymp2planwave+sphwave}}.
Combining Lemmas \ref{lm:AB_keylemma_lm1}, \ref{lm:AB_keylemma_lm3} 
and \ref{lm:AB_keylemma_lm4}, 
we get for $|x| \ge 1$
\begin{align}
\int_{\mathbb R^2} &\left\{
\frac{e^{ia|x-y|}}{|x-y|^{1/2}}-\frac{e^{ia(|x|-\omega_x\cdot y)}}{|x|^{1/2}}
\right\}u(y)dy     \nonumber\\
\noalign{\vskip 4pt}
&\leq C
\begin{cases}
|x|^{-(\sigma-1)/2}& \mbox{\rm if \ }2<\sigma < 4, \\
|x|^{-3/2}\log(1 +|x|) &  \mbox{\rm if \ }\sigma=4, \\
|x|^{-3/2}& \mbox{\rm if \ }\sigma>4, 
\end{cases} 
\end{align}
where $C$ is a positive constant independent of $a$. 
This fact, together with (\ref{eq:varphi-varphi0}) and (\ref{eq:def_f_pm}), 
gives Theorem \ref{eq:planwave+sphwave_esti}. 
\QED


\section*{Appendix}
\appendix
\section{Some inequalities}
\begin{lemma}
Let $n\in\mathbb N$ and $\Phi(x)$ be defined by
\[\Phi(x):=\int_{\mathbb R^n}\frac1{|x-y|^\beta\langle y\rangle^\gamma}dy.\]
If $0<\beta<n$ and $\beta+\gamma>n$, then $\Phi(x)$ is a bounded continuous function satisfying
\[|\Phi(x)|\leq C_{\beta\gamma n}
\begin{cases}
\langle x\rangle^{-(\beta+\gamma-n)}~~~~~&if~~0<\gamma<n,\\
\langle x\rangle^{-\beta}\log(1+\langle x\rangle)~~~~~&if~~\gamma=n,\\
\langle x\rangle^{-\beta}&if ~~ \gamma>n.
\end{cases}\]
where $C_{\beta\gamma n}$ is a constant depending on $\beta,~\gamma$ and $n$.
\label{lm:UmedaA.1}\end{lemma}

For the proof of this lemma, see \cite[Lemma A.1]{Umeda3}.

Young's inequality for convolutions is as follows (cf. \cite[P271]{Stein1}): 
\begin{lemma}
Let $h=f*g$, then
\[\|h\|_{L^q}\leq\|f\|_{L^p}\|g\|_{L^r}\]
where $1\leq p,q,r\leq\infty$ and $1/q=1/p+1/r-1$.
\label{lm:YoungInequality}
\end{lemma}

\section{Some special functions}
For the reader's convenience, we 
summarize some  properties of the Bessel function $J_0(\rho)$,
the Neumann function $N_0(\rho)$
and the Struve function $H_0(\rho)$,
whose definitions were given
 by (\ref{eq:H0}), (\ref{eq:J0}) and (\ref{eq:N0}) respectively.

\begin{lemma}
Let $\rho\in\mathbb R$. Then 
\begin{equation}
J_0(\rho)=\left(\frac2{\pi\rho}\right)^{\!1/2}
 \cos\left(\rho-\frac\pi4\right) +O(\rho^{-3/2})
\label{eq:J0inf}
\end{equation}
\begin{equation}
N_0(\rho)=\left(\frac2{\pi\rho}\right)^{\!1/2} 
\sin\left(\rho-\frac\pi4\right)+O(\rho^{-3/2})
\label{eq:N0inf}
\end{equation}
as $\rho\to \infty$.
\label{lm:J0N0inf}
\end{lemma}

Proof. By \cite[p. 199]{Watson1}, we get
\[
\begin{split}
J_0(\rho)=\left(\frac 2{\pi \rho}\right)^{1/2}
\Big[&\cos(\rho-\frac14\pi)\cdot
\Big\{(0,0)+O(\rho^{-2})\Big\}\\
&-\sin(\rho-\frac14\pi)\cdot \Big\{\frac{(0,1)}2 \rho^{-1}+O(\rho^{-3})\Big\}\Big]
\end{split}
\]
as $\rho\to \infty$, where 
\[(0,m)=\frac{\prod_{i=1}^m\{-(2i-1)^2\}}{m!\cdot2^{2m}}
=\frac{(-1)^m\{(2m-1)!!\}^2}{m!\cdot2^{2m}}.\]
Noticing $(0,0)=1, (0,1)=-1/4$, we have the asymptotic formula (\ref{eq:J0inf}).
Similarly, we have the asymptotic formula (\ref{eq:N0inf}).
\QED

\begin{lemma}
Let $\rho\in\mathbb R$. Then
\begin{equation}
{\mathbf H}_0(\rho)=\left(\frac2{\pi\rho}\right)^{\! 1/2}
 \sin\left(\rho-\frac\pi4\right) +O(\rho^{-1})
\label{eq:H0inf}
\end{equation}
as $\rho\to\infty$
\label{lm:H0inf}\end{lemma}

\vspace{5pt}

{Proof.} Noting 
\[\Gamma(k+\frac32)=\sqrt\pi\frac{(2k+1)!!}{2^{k+1}},\]
we get the following formula from the definition (\ref{eq:H0}).
\[
{\mathbf H}_0(\rho)=\frac2\pi\sum_{k=0}^\infty\frac{(-1)^k\rho^{2k+1}}{\{(2k+1)!!\}^2}
\]
Then, by \cite[p. 333]{Watson1}, we get 
\[
{\mathbf H}_0(\rho)=N_0(\rho)+\frac{(\frac12\rho)^{-1}}{\{\Gamma(1/2)\}^2}
\sum_{k=0}^{p-1}\frac{(-1)^k(\frac12)_k(2k)!}{\rho^{2k}\cdot k!}+O(\rho^{-2p-1})
\]
as $\rho\to\infty$, where, 
\[
\big(\frac12\big)_k
=\frac12\cdot\frac32\cdot\cdots\cdot\frac{2k-1}2=\frac{(2k-1)!!}{2^k}.
\]
Since 
\[\frac{(-1)^k(\frac12)_k(2k)!}{\rho^{2k}\cdot k!}=\frac{(-1)^k(2k-1)!!(2k)!}{\rho^{2k}2^k\cdot k!},\]
and 
\[(2k)!=2^kk!(2k-1)!!,\]
we get
\[
\begin{split}
{\mathbf H}_0(\rho)=&N_0(\rho)+\frac{2}{\pi}
\sum_{k=0}^{p-1}
(-1)^k \{(2k-1)!!\}^2\rho^{-2k-1}+O(\rho^{-2p-1})\\
=&N_0(\rho)+O(\rho^{-1}).
\end{split}\]
as $\rho\to\infty$. 
Finally, using Lemma \ref{lm:J0N0inf}, we obtain this lemma.
\QED


\textbf{Acknowledgements}~~
D. Wei wishes to express his sincere thanks to his family for their love. 
He also wishes to express his sincere thanks to Mr. Y. Oda for his assistance on the numerical analysis.

\vspace{10pt}

\begin{flushright}\begin{minipage}{0.5\hsize}
T.Umeda\\
Department of Mathematical Science \\
University of Hyogo\\
Shosha 2167 \\
Himeji 671-2201, Japan.~~\\
e-mail : umeda@sci.u-hyogo.ac.jp
\end{minipage}\end{flushright}

\vspace{5pt}

\begin{flushright}\begin{minipage}{0.5\hsize}
D. Wei\\
Department of Mechanical and Control Engineering\\ 
Graduate School of Science and \\Engineering \\
Tokyo Institute of Technology\\
2-12-1 S5-22 O-okayama, Meguro-ku, Tokyo 152-8550, Japan\\
e-mail : dabi@ok.ctrl.titech.ac.jp
\end{minipage}\end{flushright}


\begin{thebibliography}{88}
\addcontentsline{toc}{section}{References}

\bibitem{Agmon1}
S. Agmon,  
\textit {Spectral properties of Schr\"odinger operators and scattering theory}, 
Ann. Scoula Norm. Sup. Pisa \textbf{4-2} (1975), 151-218.


\bibitem{Ben1}
M. Ben-Artzi and J. Nemirovski,  
\textit {Remarks on relativistic Schr\"odinger operators and their extensions}, 
Ann. Inst. Henri Poincar\'e, Phys. th\'eor. \textbf{67} (1997), 29-39.

\bibitem{Enss1}
V. Enss,  
\textit {Asymptotic completeness for quantum-mechanical potential scattering I, Short range potentials}, 
Commum. Math. Phys. \textbf{61} (1978), 258-291.

\bibitem{Erdelyi}
A. Erd\'elyi,
\textit{Tables of integral transforms}, Vol. I, McGraw-Hill (1954).

\bibitem{Ikebe2}
T. Ikebe, 
\textit {Eigenfunction expansions associated with 
the Schr\"odinger operators and their applications to scattering theory}, 
Arch. Rational Mech. Anal. \textbf{5} (1960), 1-34.

\bibitem{Isozaki1}
H. Isozaki, 
\textit {Many-body Schr\"odinger equation}, 
Springer Tokyo (2004). (In Japanese).

\bibitem{KatoKuroda}
T. Kato and S.T. Kuroda,
\textit {Theory of simple scattering and eigenfunction expansions}, 
F.E. Browder ed., Functional Analysis and Related Rields,
 Springer (1970), 99-131. 

\bibitem{Kitada3}
A. Jensen and H. Kitada, 
\textit {Fundamental solutions and eigenfunction expansions for 
Schr\"odinger Operators II: Eigenfunction Expansions}, 
Math. Z. \textbf{199} (1988), 1-13. 

\bibitem{Kuroda1}
S.T. Kuroda,  
\textit {Spectral theory II}, 
Iwanami Shoten (1979). (In Japanese).

\bibitem{Moriguti1}
S. Moriguti, K. Udagawa and S. Hitotsumatu
\textit {Mathematical fomula II}, 19th Edition, 
Iwanami Syoten (2002). (In Japanese).

\bibitem{Moriguti2}
S. Moriguti, K. Udagawa and S. Hitotsumatu
\textit {Mathematical formula III}, 19th Edition, 
Iwanami Syoten (2002) (In Japanese).

\bibitem{Reed1}
M. Reed and B. Simon, 
\textit {Methods of Modern Mathematical Physics IV: Analysis of Operators}, 
Academic Press (1978).

\bibitem{Simon1}
B. Simon, 
\textit {Phase space analysis of simple scattering systems: Extensions of some work of Enss}, 
Duke Math. J. \textbf{46} (1979), 119-168.

\bibitem{Stein1}
E.M. Stein, 
\textit {Singular integrals and differentiability properties of functions}, 
Princeton University Press (1970).

\bibitem{Umeda2}
T.Umeda, 
\textit {The action of $\sqrt{-\Delta}$ on weighted Sobolev spaces}, 
Lett. Math. Phys. {\bf 54} (2000), 301-313.

\bibitem{Umeda5}
T.Umeda, 
\textit {Eigenfunction expansions associated with relativistic Schr\"odinger operators}, 
in Partial Differential Equations and Spectral Theory, eds. M. Demuth and B. W. Schulze, Operator Theory: 
Advances and Applications {\bf 126} (2001), 315-319.


\bibitem{Umeda3}
T. Umeda, 
\textit {Generalized eigenfunctions of relativistic Schr$\ddot{o}$dinger operators I}, 
Electron. J. Diff. Eqns. {\bf 2006} (2006), 1-46.

\bibitem{Watson1}
G. N. Watson, 
\textit{Theory of Bessel Functions}, 2nd Edition, 
Cambridge  University Press  (1966).

\bibitem{Wei1}
D. Wei, 
\textit {Completeness of the Generalized Eigenfunctions for relativistic Schr\"odinger operators I}, 
Osaka J. Math. {\bf 44} (2007), 851-881.

\bibitem{Wei3}
D. Wei, 
\textit {Completeness of the Generalized Eigenfunctions for relativistic Schr\"odinger operators II}, 
in preparation.

\end{thebibliography}
\end{document}